\documentclass[12pt]{article}
\usepackage[english]{babel}
\usepackage{amsmath,cite,amsthm}
\usepackage{amssymb}
\usepackage{latexsym}
\newtheoremstyle{theorem}
  {10pt}          
  {10pt}  
  {\sl}  
  {\parindent}     
  {\bf}  
  {. }    
  { }    
  {}     
\theoremstyle{theorem}
\newtheorem{theorem}{Theorem}

\newtheoremstyle{defi}
  {10pt}          
  {10pt}  
  {\rm}  
  {\parindent}     
  {\bf}  
  {. }    
  { }    
  {}     
\theoremstyle{defi}

\begin{document}

\title{$\mathbb{B}$-valued monogenic functions 
and their applications to boundary value problems in displacements
of\\ 2-D Elasticity}

\author{S.V. Gryshchuk\\
 Institute of Mathematics, \\ National Academy of Sciences
of Ukraine,\\ Tereshchenkivska Str. 3, 01601, Kyiv, Ukraine\\
gryshchuk@imath.kiev.ua, serhii.gryshchuk@gmail.com}

\maketitle

\begin{abstract}
Consider the commutative algebra $\mathbb{B}$ over the field of
complex numbers with the bases $\{e_1,e_2\}$ such that  
$(e_1^2+e_2^2)^2=0$, $e_1^2+e_2^2\ne 0$. 
 Let $D$ be a domain in $xOy$,
$D_{\zeta}:=\{xe_1+ye_2:(x,y) \in  D\}\subset \mathbb{B}$. We say
that $\mathbb{B}$-valued function  $\Phi \colon D_{\zeta}
\longrightarrow \mathbb{B}$, $\Phi(\zeta)=U_{1}\,e_1+U_{2}\,ie_1+
U_{3}\,e_2+U_{4}\,ie_2$, $\zeta=xe_1+ye_2$, $U_{k}=U_{k}(x,y)\colon
D\longrightarrow \mathbb{R}$, $k=\overline{1,4}$, is {\em monogenic}
in $D_{\zeta}$ iff $\Phi$ has the classic derivative in every point
in $D_{\zeta}$. Every $U_k$, $k=\overline{1,4}$, is a biharmonic
function in $D$. A problem on finding an elastic equilibrium for
isotropic body $D$ by given boundary values on $\partial D$ of
partial derivatives $\frac{\partial u}{\partial v}$, $\frac{\partial
v}{\partial y}$ for displacements $u$, $v$ is equivalent to BVP for
monogenic functions, which is to find  $\Phi$ by given boundary
values of $U_1$ and $U_4$.

 \vspace*{6mm}


\end{abstract}

\newpage

\section{Monogenic functions in the bi\-har\-mo\-nic algebra asso\-cia\-ted with the bi\-har\-mo\-nic equ\-ati\-on}                

We say that an associative commutative two-dimensional algebra
$\mathbb B$ with the unit $1$ over the field of complex numbers
$\mathbb C$ is {\it biharmonic} if in $\mathbb B$ there exists a
{\it biharmonic} basis, i.e, a bases $\{e_1,e_2\}$ satisfying the
conditions
\begin{equation}\label{biharm-bas}
 (e_1^2+e_2^2)^2=0,\qquad e_1^2+e_2^2\ne 0\,.
\end{equation}

V.\,F. Kovalev and I.\,P. Mel'nichenko \cite{Kov-Mel}
 found a
multiplication table for a biharmonic basis $\{e_1,e_2\}$:
\begin{equation} \label{tab_umn_bb}
e_1=1,\qquad e_2^2=e_1+2ie_2,
\end{equation}
where $i$\,\, is the imaginary complex unit.

In \cite{Mel86}, I.~P.~Mel'nichenko proved that there exists the
unique biharmonic algebra $\mathbb{B}$ and he constructed
all biharmonic bases. 

Consider a {\it biharmonic plane} $\mu:=\{\zeta=x\,e_1+y\,e_2 :
x,y\in\mathbb R\}$ which is a linear span of the elements $e_1,e_2$
of the biharmonic basis (\ref{tab_umn_bb}) over the field of real
numbers $\mathbb R$. With a domain $D$ of the Cartesian plane $xOy$
we associate the congruent domain $D_{\zeta}:= \{\zeta=xe_1+ye_2 :
(x,y)\in D\}$ in the biharmonic plane $\mu$, and corresponding
domain in the complex plane $\mathbb{C}$: $D_{z}:=\{z=x+iy:(x,y)\in
D\}$.

Let $D_{\ast}$ be a domain in $xOy$ or in $\mu$.
 Denote by $\partial D_{\ast}$
 a boundary of a
domain $D_{\ast}$, $\mathrm{cl}D_{\ast}$ means a closure of a domain
$D_{\ast}$.

In what follows, $(x,y)\in D$, $\zeta=x\,e_1+y\,e_2\in D_{\zeta}$,
$z=x+iy\in D_{z}$.

Inasmuch as divisors of zero don't belong to the biharmonic plane,
one can define the derivative $\Phi'(\zeta)$ of function $\Phi
\colon D_{\zeta}\longrightarrow \mathbb{B}$ in the same way as in
the complex plane:
\[\Phi'(\zeta):=\lim\limits_{h\to 0,\, h\in\mu}
\bigl(\Phi(\zeta+h)-\Phi(\zeta)\bigr)\,h^{-1}\,.\]
 We say that a function $\Phi : D_{\zeta}\longrightarrow \mathbb{B}$ is
\textit{monogenic} in a domain $D_{\zeta}$ and, denote by
$\Phi\in\mathcal{M}_{\mathbb B}(D_{\zeta})$, iff the derivative
$\Phi'(\zeta)$ exists in every point $\zeta\in D_{\zeta}$.


Every function  $\Phi \colon D_{\zeta}\longrightarrow \mathbb{B}$
has a form
\begin{equation}\label{mon-funk} \Phi(\zeta)=U_{1}(x,y)\,e_1+U_{2}(x,y)\,ie_1+
U_{3}(x,y)\,e_2+U_{4}(x,y)\,ie_2,
\end{equation}
where $\zeta=xe_1+ye_2$, $U_{k}\colon D\longrightarrow \mathbb{R}$,
$k=\overline{1,4}$.

Every real component  $U_{k}$, $k=\overline{1,4}$, in expansion
 (\ref{mon-funk}) we denote by
$\mathrm{U}_{k}\left[\Phi\right]$, i.e., for 
$k\in\{1,\dots,4\}$:
$\mathrm{U}_{k}\left[\Phi(\zeta)\right]:=U_{k}(x,y)$ for all $\zeta=
xe_1 + y e_2 \in D_{\zeta}$.


It is established in the paper   \cite{Kov-Mel} that a function
$\Phi\colon D_{\zeta}\longrightarrow \mathbb{B}$ is monogenic in a
domain $D_{\zeta}$ if and only if components $U_{k}$,
$k=\overline{1,4}$, of the expression \eqref{mon-funk} are
differentiable in the domain $D$ and the following analog of the
Cauchy -- Riemann conditions  is satisfied:
\begin{equation}\label{usl_K_R}
\frac{\partial \Phi(\zeta)}{\partial y}\,e_1=\frac{\partial
\Phi(\zeta)}{\partial x}\,e_2.
\end{equation}

In an extended form the condition (\ref{usl_K_R}) for the monogenic
function (\ref{mon-funk}) is equivalent to the system of four
equations (cf., e.g., \cite{Kov-Mel,GrPl_umz-09}) with respect to
components $U_{k}=\mathrm{U}_{k}\left[\Phi\right]$,
$k=\overline{1,4}$, in (\ref{mon-funk}):
\begin{eqnarray}
\frac{\partial U_{1}(x,y)}{\partial y}&=&\frac{\partial
U_{3}(x,y)}{\partial x},\label{kr1}\\ 
\frac{\partial U_2(x,y)}{\partial y}&=&\frac{\partial U_{4}(x,y)}{\partial x},\label{kr2}\\
 \frac{\partial U_{3}(x,y)}{\partial y}&=&\frac{\partial U_{1}(x,y)}{\partial x}-2\frac{\partial U_{4}(x,y)}{\partial x},\label{kr3}\\
\frac{\partial U_{4}(x,y)}{\partial y}&=&\frac{\partial
U_{2}(x,y)}{\partial x}+2\frac{\partial U_{3}(x,y)}{\partial
x}\label{kr4}.
\end{eqnarray}

It is proved in the paper \cite{Kov-Mel} that a function
$\Phi(\zeta)$ having derivatives till fourth order in $D_{\zeta}$
satisfies the two-dimensional biharmonic equation
\begin{equation}\label{big-eq}
(\Delta_2)^{2}U(x,y):= \left(\frac{\partial^4}{\partial
x^4}+2\,\frac{\partial^4}{\partial x^2\partial
y^2}+\frac{\partial^4}{\partial y^4}\right)U(x,y)=0
\end{equation}
in the domain $D$ owing to the relations (\ref{biharm-bas}) and
\[(\Delta_2)^{2}\Phi(\zeta)=\Phi^{(4)}(\zeta)\,(e_1^2+e_2^2)^2.\]
Therefore, every component $U_{k}\colon D\longrightarrow
\mathbb{R}$, $k=\overline{1,4}$, of the expansion  (\ref{mon-funk})
satisfies also the equation (\ref{big-eq}), i.~e. $U_{k}$ is a {\em
biharmonic function} in the domain $D$.

It is proved  \cite{GrPl_umz-09} that a monogenic function $\Phi
\colon D_{\zeta}\longrightarrow \mathbb{B}$ has derivatives
$\Phi^{(n)}(\zeta)$ of all orders in the domain $D_{\zeta}$ and,
consequently, satisfies the two-dimensional biharmonic equation
(\ref{big-eq}).

 In the papers
\cite{GrPl_umz-09,Cont_2011,Cont-13,Pl_Brig} it was also proved such
a fact that every biharmomic function $U_{1}(x,y)$ in a bounded
simply connected domain $D$  is the first component of the expansion
(\ref{mon-funk}) of monogenic function $\Phi \colon
D_{\zeta}\longrightarrow \mathbb{B}$,
moreover, all such functions $\Phi$ are found 
in an explicit form.

Basic analytic properties of monogenic functions in a biharmonic
plane are similar to properties of holomorphic functions of the
complex variable. More exactly, analogues of the Cauchy integral
theorem and integral formula, the Morera theorem, the uniqueness
theorem, the Taylor and Laurent expansions are established in papers
\cite{Cont-13,Pl_Brig,Gr-Pl_Dop2009}.

\section{Auxiliary statement}
\begin{theorem}  \label{2der-bih-f+mon-f}
Let $W$ be an arbitrary fixed biharmonic in a domain  $D$ function;
${\Phi}_{\ast}\in\mathcal{M}_{\mathbb B}(D_{\zeta})$,
${\Phi}\in\mathcal{M}_{\mathbb B}(D_{\zeta})$ and
$\mathrm{U}_{1}\left[{\Phi}_{\ast}\right]=W$,
$\Phi:={\Phi}_{\ast}''$.
 Then the following formulas are true:
\begin{equation}\label{form-2-der}
\frac{\partial^{2}W(x,y)}{\partial {x}^{2}}=
\mathrm{U}_{1}\left[\Phi(\zeta)\right], \,
  \frac{\partial^{2}W(x,y)}{\partial {y}^{2}}= \mathrm{U}_{1}\left[\Phi(\zeta)\right]-2
  \mathrm{U}_{4}\left[\Phi(\zeta)\right], 
\end{equation}
for every $(x,y)\in D$, $\zeta=xe_1+ye_2\in D_{\zeta}$.
\end{theorem}

\vspace{3mm} \noindent \textbf{Proof.} There exists
${\Phi}_{\ast}\in\mathcal{M}_{\mathbb B}(D_{\zeta})$ such that
\begin{equation}\label{1-Phi_ast}
\mathrm{U}_{1}\left[{\Phi}_{\ast}(\zeta)\right]=W(x,y)\qquad \forall
\zeta \in D_{\zeta}.
\end{equation}

A derivative of a function  $\Phi_{\ast}$ is represented by the
equality  $\Phi_{\ast}'=\frac{\partial \Phi_{\ast}}{\partial x}$.
Therefore, we have the following equalities
\begin{equation}\label{comp-Phi_ast}
\mathrm{U}_{k}\left[{\Phi}_{\ast}'(\zeta)\right]=\frac{\partial
\mathrm{U}_{k}\left[\Phi_{\ast}(\zeta)\right]}{\partial x}, k \in
\{1,\dots,4\} , \qquad \forall \zeta \in D_{\zeta}.
\end{equation}

Using the equality  (\ref{kr1}) for monogenic function
$\Phi_{\ast}$, deliver the equality $$\frac{\partial
\mathrm{U}_{1}\left[{\Phi}_{\ast}(\zeta)\right]}{\partial y}=
\frac{\partial
\mathrm{U}_{3}\left[{\Phi}_{\ast}(\zeta)\right]}{\partial x}\qquad
\forall \zeta \in D_{\zeta},$$ substituting into which successively
equalities
 (\ref{comp-Phi_ast}) with
$k=3$ and (\ref{1-Phi_ast}), as a result, obtain
\begin{equation}\label{comp-3-Phi_ast'}
\mathrm{U}_{3}\left[{\Phi}_{\ast}'(\zeta)\right]=\frac{\partial
W(x,y)}{\partial y} \qquad \forall \zeta \in D_{\zeta}.
\end{equation}

Now, substituting  (\ref{1-Phi_ast}) in (\ref{comp-Phi_ast}) with
$k=1$, obtain
\begin{equation}\label{comp-1-Phi_ast'}
\mathrm{U}_{1}\left[{\Phi}_{\ast}'(\zeta)\right]=\frac{\partial
W(x,y)}{\partial x} \qquad \forall \zeta \in D_{\zeta}.
\end{equation}

Since, $\Phi \equiv \Phi_{\ast}''=\frac{\partial
\Phi_{\ast}'}{\partial x}$, therefore, it implies the equalities
\begin{equation}\label{comp-Phi_ast'}
\mathrm{U}_{k}\left[\Phi (\zeta)\right]=\frac{\partial
\mathrm{U}_{k}\left[\Phi_{\ast}'(\zeta)\right]}{\partial x}, k \in
\{1,\dots,4\} , \qquad \forall \zeta \in D_{\zeta}.
\end{equation}

Substituting  (\ref{comp-1-Phi_ast'}) in the equality
(\ref{comp-Phi_ast'}) with  $k=1$, we get 
the first equality in (\ref{form-2-der}).

Using the formula  (\ref{kr3}) for monogenic function
$\Phi_{\ast}'$, we have
\begin{equation}\label{sootn}
 \frac{\partial
\mathrm{U}_{3}\left[\Phi_{\ast}'(\zeta)\right]}{\partial
y}=\frac{\partial
\mathrm{U}_{1}\left[\Phi_{\ast}'(\zeta)\right]}{\partial
x}-2\frac{\partial
\mathrm{U}_{4}\left[\Phi_{\ast}'(\zeta)\right]}{\partial x} \qquad
\forall \zeta \in D_{\zeta}. \end{equation}

 Finally, substituting in series 
  relations
(\ref{comp-3-Phi_ast'}), (\ref{comp-Phi_ast'}) with $k=1$ and $k=4$
into (\ref{sootn}), obtain the second equality in
(\ref{form-2-der}). The theorem is proved.

\section{Displacements-type problem}
We shall assume further that  $D$ is a bounded simply connected
domain in the Cartesian plane $xOy$. 
Let $\tau \colon D\longrightarrow \mathbb{R}$ be a real-valued
function, then by ${\mathcal
C}^{k}(D)$, $k \in \{0, 1,\dots\}$, we denote a class of 
functions having continuous derivatives up to the
 $k$-th order inclusively, 
 ${\mathcal C}(D):= {\mathcal C}^{0}(D)$. If $\tau\in {\mathcal C}^{2}(D)$ then $\Delta_2
\tau:=\frac{\partial^{2}\tau}{\partial
x^2}+\frac{\partial^{2}\tau}{\partial y^2}$. For a function $\tau
\in {\mathcal C}\left(\mathrm{D}\right)$ denote by $\tau \in
{\mathcal C}\left(\mathrm{cl}{D}\right)$ if there exists a finite
limit $$\left.\tau(x,y)\right|_{(x_0,y_0)}:= \lim\limits_{(x,y)\in
D, (x,y) \to (x_0,y_0)} \tau(x,y)\qquad \forall (x_0,y_0) \in
\partial D.$$


Let a function $\Phi$ is of the type  $\Phi \colon D_{\zeta}
\longrightarrow \mathbb{B}$. Denote by $\Phi \in {\mathcal
C}\left(\mathrm{cl}{D_{\zeta}}\right)$
 if and only if $\mathrm{U}_{k}\left[\Phi\right]\in
{\mathcal C}\left(\mathrm{cl}D\right)$, $k \in
\{1,\dots,4\}$. 
 In this case we use for every $\zeta_0=x_0 e_1+ y_0 e_2 \in\partial
D_{\zeta}$  notations
$$\mathrm{U}_{k}\left[\Phi(\zeta_0)\right]:=
\mathrm{U}_{k}\left[\lim\limits_{\zeta \in D_{\zeta}, \zeta \to
\zeta_0\in\partial D_{\zeta}}{\Phi}(\zeta)\right], k
=\overline{1,4}.$$

Consider a  boundary value problem: to find  in $D$ partial
derivatives $\mathcal{V}_{1}:=\frac{\partial u}{\partial x}$,
$\mathcal{V}_{2}:=\frac{\partial v}{\partial y}$ for displacements
 $u=u(x,y)$, $v=v(x,y)$ of an elastic isotropic body occupying $D$,
when   their boundary values are  given on the boundary $\partial
D$: 
\begin{equation}\label{u_x-v_y}
\left.\mathcal{V}_{k}(x,y)\right|_{(x_0, y_0)}= g_{k}(x_0, y_0),
k=1,2, \quad \forall (x_0, y_0)\in \partial D,
\end{equation}
where $g_{k} \colon \partial D \longrightarrow \mathbb{R}$, $k=1,
2$, are given functions. 


We shall call this problem as the  {\em  $(u_{x},v_{y})$-problem}.
Let $W \colon D \longrightarrow \mathbb{R}$ be an unknown biharmonic
function. Further we mean by this function the Airy stress function. 
Denote 
\begin{equation}\label{C_1-2}
 \mathrm{C}_{k}[W](x,y):=  W_{k}(x,y) +\kappa_{0} \,W_{0}(x,y) , k=1,2, \qquad  \forall (x,y)\in\mathrm{cl}{D},
\end{equation}
where
\begin{equation}\label{W_0}
W_{0}(x,y):= \Delta_2 W(x,y),
\end{equation}
\begin{equation}\label{W_1-2} W_{1}(x,y):=\frac{\partial^2
W(x,y)}{\partial x^2},\, W_{2}(x,y):=\frac{\partial^2
W(x,y)}{\partial y^2},
\end{equation}
$\kappa_{0}:= \frac{\lambda+2\mu}{2(\lambda+\mu)}$ and
 $\mu$, $\lambda$ are
Lam\'{e} constants (cf., e.g., \cite[p.~2]{Lu}).

The following equalities are valid in $D$  
(cf., e.g., \cite[pp.~8~--~9]{Lu},\cite[p.~5]{Mikhlin}):
\begin{equation}\label{Syst-14}
2\mu\,\mathcal{V}_{k}(x,y)=\mathrm{C}_{k}[W](x, y), k=1,2,  \qquad
\forall (x, y)\in D.
\end{equation}

 Then solving the $(u_{x},v_{y})$-problem is equivalent 
 to finding
in $D$ quantities of  $\mathrm{C}_{k}[W]$, $k=1,2$, where boundary
values of an unknown biharmonic function  $W \colon D
\longrightarrow \mathbb{R}$ satisfy the system
\begin{equation}
\label{Syst-14-BVP} 2\mu\, \mathrm{C}_{k}[W](x_0, y_0)= g_{k}(x_0,
y_0), k=1,2, \qquad \forall (x_0, y_0)\in
\partial D.
\end{equation}

\begin{theorem}\label{2-part-W-BVP}
The $(u_{x},v_{y})$-problem is equivalent to boundary value problem
on finding in  $D$ the second derivatives $\frac{\partial^{2}
W(x,y)}{\partial x^2}$, $\frac{\partial^{2} W(x,y)}{\partial y^2}$
of a biharmonic function  $W \in {\mathcal C}^{2}(\mathrm{cl}{D})$,
which satisfy for all  $(x_0,y_0)\in
\partial D$ the boundary data\emph{:}
\begin{equation}
\label{bv-W-1} \left.\frac{\partial^{2} W(x,y)}{\partial
x^2}\right|_{(x_0,y_0)}= \lambda\,
g_{1}(x_0,y_0)+(\lambda+2\mu)\,g_{2}(x_0,y_0),
\end{equation}
\begin{equation}
\label{bv-W-2} \left.\frac{\partial^{2} W(x,y)}{\partial
y^2}\right|_{(x_0,y_0)}= (\lambda+2\mu)\,
g_{1}(x_0,y_0)+\lambda\,g_{2}(x_0,y_0).
\end{equation}
Then a general solution of the $(u_{x},v_{y})$-problem is expressed
by the formulas\emph{:}
\begin{equation}\label{Sot-W-2}
 \mathcal{V}_{k}(x,y)=\frac{1}{2\mu}\, \mathrm{C}_{k}[W](x,y), k=1,2, \qquad \forall (x,y) \in D.
 \end{equation}
\end{theorem}

\vspace{3mm} \noindent \textbf{Proof.}  Adding equalities of the
system
 (\ref{Syst-14}) term by term,
taking into account definitions (\ref{C_1-2}) and the value of  $\kappa_{0}$, 
we get
\begin{equation}\label{Delta-W}
2\mu \left(\mathcal{V}_{1}(x,y)+\mathcal{V}_{2}(x,y)\right)=
-\frac{\lambda}{2(\lambda+\mu)}\Delta_2 W(x,y)\qquad \forall
(x,y)\in D.
\end{equation}
Now, it follows from formulas  (\ref{Syst-14}) and (\ref{Delta-W}),
definitions
(\ref{C_1-2}),  
 inclusions 
 $\mathcal{V}_{k}\in
\mathcal{C}(\mathrm{cl}D)$, $k=1,2$, 
that functions $\Delta_2
W$, $\frac{\partial^2 W}{\partial x^2}$, $\frac{\partial^2
W}{\partial y^2}$ belong to $\mathcal{C}\left(\mathrm{cl}D\right)$,
therefore, $W\in\mathcal{C}^{2}\left(\mathrm{cl}D\right)$.

Now, solving the  system (\ref{Syst-14-BVP}) with respect to
$\left.\frac{\partial^{2} W(x,y)}{\partial x^2}\right|_{(x_0,y_0)}$
and $\left.\frac{\partial^{2} W(x,y)}{\partial
y^2}\right|_{(x_0,y_0)}$, we obtain the system of equations
(\ref{bv-W-1}),(\ref{bv-W-2}). Formulas  (\ref{Sot-W-2}) are
followed the equalities (\ref{Syst-14}).

In a similar way, we can prove that solving of a boundary value
problem for biharmonic function $W$ with boundary data
(\ref{bv-W-1}), (\ref{bv-W-1}) implies a solution of the
$(u_{x},v_{y})$-problem by the formulas (\ref{Sot-W-2}), where a
function (\ref{W_0}) is a sum of functions in (\ref{W_1-2}).
 The theorem is proved.

\begin{theorem}\label{hom-u_x-v_y}
 A general
solution of the homogeneous $(u_{x},v_{y})$-problem with zero data
$g_1=g_2\equiv 0$ is the trivial\emph{:}
\begin{equation}\label{sol-odn-ux-vy}
\mathcal{V}_{k}(x,y)\equiv 0, k=1,2, \qquad \forall (x,y) \in D.
\end{equation}
\end{theorem}

\vspace{3mm} \noindent \textbf{Proof.} 
 Passing in (\ref{Delta-W}) to
a limit, as  $(x,y)$ tends to arbitrary fixed point
$(x_0,y_0)\in\partial D$, and taking
into account definitions  (\ref{W_0}), we conclude that  
$\left.W_{0}(x,y)\right|_{(x_0,y_0)}=0$.
 Inasmuch as, (\ref{W_0}) is a harmonic function then from the last boundary
 equality, we obtain by the maximum principle the equality
\begin{equation}\label{widetilde-W}
W_{0}(x,y)= 0 \qquad \forall (x,y)\in\mathrm{cl}{D}.
\end{equation}

The equality (\ref{widetilde-W}) yields that functions (\ref{W_1-2})
 are harmonic in  $D$ and, by formulas (\ref{bv-W-1}) and
(\ref{bv-W-2}), vanish on the boundary $\partial D$, thus, in view
of the maximum principle, we have that
\begin{equation}\label{W-1_W-2}
W_1(x,y)=W_2(x,y)\equiv 0\quad \forall (x,y) \in D.
 \end{equation}
 Substituting in series  expressions  (\ref{C_1-2}), 
equalities  (\ref{widetilde-W}) and (\ref{W-1_W-2}) into the formulas  (\ref{Sot-W-2}), we get the 
formulas  (\ref{sol-odn-ux-vy}). The theorem is proved.

Note, that in \cite{KM-88-sys-Lame,Gr-dop-15} considered expressions
of solutions the Lam\'{e} equilibrium system in displacements via
components $U_{k}$, $k=\overline{1,4}$, of a monogenic function
(\ref{mon-funk}). Moreover, the following statemant is proved in
\cite[Theorem~1]{Gr-dop-15}:

\emph{Let a function  \emph{(\ref{mon-funk})} is monogenic in a
domain
$D_{\zeta}$. 
Then the next pairs of functions
\[u(x,y)=\frac{2}{\gamma}\,U_{1}(x,y)-\frac{2+\gamma}{\gamma}\,U_{4}(x,y), \,  v(x,y)=U_{2}(x,y); \]
\[u(x,y)=-\frac{2+\gamma}{\gamma}\,U_{2}(x,y)-\frac{2(1+\gamma)}{\gamma}\,U_{3}(x,y), \,  v(x,y)=U_{4}(x,y);\]
\[u(x,y)=-\frac{2}{\gamma}\,U_{2}(x,y)-\frac{2+\gamma}{\gamma}\,U_{3}(x,y), \,  v(x,y)=U_{1}(x,y)
\]
are solutions the Lam\'{e} equilibrium system in displacements
\begin{equation}
    \label{Lame}     
\left\{
\begin{array}{ll}
\Delta u +\gamma \frac{\partial \theta}{\partial x}=0,\vspace*{2mm}\\
\Delta v +\gamma \frac{\partial \theta}{\partial y}=0,
\end{array}
\right.
\end{equation}
where $\theta:=\frac{\partial u}{\partial x}+\frac{\partial
v}{\partial y}$,  $\gamma:=(\lambda+\mu)\mu^{-1}$.}  

Note, that this statement is generalization of a result in
\cite{KM-88-sys-Lame} to a general bounded 
domain $D_{\zeta}$. Another results of \cite{KM-88-sys-Lame} are
also generalized to this case. 
\section{Boundary value (1-4)-problem for monogenic functions}
V.\,F.~Kovalev \cite{Kov,IJPAM_13} posed the {\it biharmonic Schwarz
type problems} on finding\break $\Phi\in\mathcal{M}_{\mathbb
B}(D_{\zeta})\cap
{\mathcal C}(\mathrm{cl}D_{\zeta})$ by given boundary values of $U_k$ and $U_m$, $1\le k<m\le 4$, in (\ref{mon-funk}) : 
\[U_{k}(x,y)=u_{k}(\zeta)\,,\quad
U_{m}(x,y)=u_{m}(\zeta)\qquad\forall\, \zeta \in\partial
D_{\zeta},\,\] where $u_{k}$ and $u_{m}$ are given real-valued
functions. we We shall call this problem by the $(k-m)$-problem

Some relations between the (1-3)-problem and problems of the theory
of elasticity are described in
\cite{Kov,G_conf_14,G-P_arXiv-15,G-P_MMAS-mf_BVP}. In particular, it
is shown that the main biharmonic problem (cf., e.g., \cite[p.
194]{S3} and \cite[p. 13]{Mikhlin}) on finding a biharmonic function
$U \colon  D \longrightarrow\mathbb{R}$ with given limiting values
of its partial derivatives $\left.\frac{\partial U}{\partial
x}\right|_{(x_0,y_0)}$ and  $\left.\frac{\partial U}{\partial
y}\right|_{(x_0,y_0)}$  can be reduced to the (1-3)-problem.

In  \cite{IJPAM_13,G_P-Proc-Mosc-11}, we investigated the
(1-3)-problem for cases where $D_{\zeta}$ is either an upper
half-plane or a unit disk in the biharmonic plane. Its solutions
were found in explicit forms with using of some integrals analogous
to the classic Schwarz integral. Moreover, the (1-3)-problem is
solvable unconditionally for a half-plane but it is solvable for a
disk if and only if a certain condition is satisfied.

In \cite{G_conf_14}, a certain scheme was proposed for reducing the
(1-3)-problem in a simply connected domain with sufficiently smooth
boundary to a suitable boundary value problem in a disk with using
power series and conformal mappings in the complex plane.

Under general suitable smooth conditions for a boundary of a bounded
domain $D$ and using a hypercomplex analog of the Cauchy type
integral, we reduce the (1-3) boundary value problem to a system of
integral equations on the real axes and establish sufficient
conditions under which this system has the Fredholm property
\cite{G-P_arXiv-15,G-P_MMAS-mf_BVP}.

In this section we are interested in the (1-4)-problem with boundary
data
\begin{equation} \label{u_1-4}
\mathrm{U}_{k}\left[\Phi(\zeta_0)\right]=u_{k}(\zeta_0), k \in
\{1,4\},\quad \forall\, \zeta_0=x_0 e_1+ y_0 e_2 \in \partial
D_{\zeta}.
\end{equation}
\begin{theorem}\label{syst-1-4} Let $W$ is  a general biharmonic function from  ${\mathcal
C}^{2}(\mathrm{cl}D)$ satisfying the boundary conditions
\emph{(\ref{Syst-14-BVP})}. Then $W$ rebuilds
 a general solution of the  \emph{($u_x, v_y$)}-problem with boundary data
\emph{(\ref{u_x-v_y})} by the formulas \emph{(\ref{Syst-14})}.

 A general solution $\Phi$ of the (1-4)-problem with boundary data    
 \begin{equation}\label{u_1-4}
 u_1=\lambda\, g_1+\left(\lambda+2\mu\right)\,g_2,
 u_4=-\mu\,g_1+\mu\, g_2,
 \end{equation}
 generate in  $D$ the second order derivatives  $\frac{\partial^2 W}{\partial
 x^2}$, $\frac{\partial^2 W}{\partial
 y^2}$ 
 by the formulas
 \emph{(\ref{form-2-der})}. A general solution of the $(u_{x},v_{y})$-problem for any
 $(x,y)\in \mathrm{cl}{D}$
 is expressed by the equalities
 \begin{equation}\label{u-x}
 2\mu \frac{\partial u(x,y)}{\partial x}=
 \frac{\mu}{\lambda+\mu}\mathrm{U_1}\,\left[\Phi(\zeta)\right]
 -\frac{\lambda+2\mu}{\lambda+\mu} \mathrm{U_4}\,\left[\Phi(\zeta)\right],
 \end{equation}
 \begin{equation} \label{v-y}
2\mu \frac{\partial v(x,y)}{\partial y}=
\frac{\mu}{\lambda+\mu}\mathrm{U_1}\,\left[\Phi(\zeta)\right]+
 \frac{\lambda+2\mu}{\lambda+\mu}\mathrm{U_4}\,\left[\Phi(\zeta)\right],
  \end{equation}
where $\zeta=x e_1+ y e_2 \in \mathrm{cl}{D_{\zeta}}$.
\end{theorem}

\vspace{3mm} \noindent \textbf{Proof.} By Theorem~
\ref{2-part-W-BVP}, the $(u_{x},v_{y})$-problem is equivalent to
finding the second order partial derivatives $\frac{\partial^2
W}{\partial x^2}$, $\frac{\partial^2 W}{\partial y^2}$ of a
sought-for biharmonic function\break  $W\in{\mathcal
C}^{2}(\mathrm{cl}{D})$,  which satisfy the limiting conditions
(\ref{bv-W-1}) and (\ref{bv-W-2}).

Let $\Phi$ be a required general solution of the  (1-4)-problem with boundary data   
 (\ref{u_1-4}). Then, taking into account that  $D$ is a simply-connected,
 obtain that there exists a function  $\Phi_{\ast}\in\mathcal{M}_{\mathbb
 B}(D_{\zeta})$ such that
 $\Phi(\zeta)=\Phi_{\ast}''(\zeta)$ for all $\zeta\in D_{\zeta}$
 and $\mathrm{U}_{1}[\Phi_{\ast}]= W$. 
  Then by Theorem~\ref{2der-bih-f+mon-f}  the relations (\ref{form-2-der}) are true,
substitute them into the expressions  (\ref{C_1-2}). 

After computations, we obtain that $\mathrm{C}_{1}[W]$ is equal to
the right-part in  (\ref{u-x}), $\mathrm{C}_{2}[W]$ is equal to the
right-part in (\ref{v-y}). Passing to a limit in these equalities,
as $(x,y)$ tends to arbitrary boundary point  $(x_0, y_0)\in
\partial D_{\zeta}$, we get with use of equalities
(\ref{Syst-14-BVP}), (\ref{u_1-4}) an equatily.

Substituting into the right-part in (\ref{Syst-14}) formulas
(\ref{form-2-der}), we obtain the equalities (\ref{u-x}),
(\ref{v-y}). The theorem is proved.


\section{$(u_{x},v_{y})$-problem and  the elastic
equilibrium}

We want to find how a solution of the $(u_{x},v_{y})$-problem
generates stresses $\sigma_x$, $\tau_{xy}$, $\sigma_y$. Assume that
we know a general solution of the $(u_{x},v_{y})$-problem\break
 $\mathcal{V}_{1}=\frac{\partial u}{\partial x}$,
$\mathcal{V}_{2}=\frac{\partial v}{\partial y}$. Then an unknown the
Airy stress function $W$ need to satisfy
 conditions (\ref{Syst-14}).
 By the generalized Hooke's law, we have the system of three equations:
\begin{equation}\label{sig-x}
 \sigma_x =  (\lambda+2\mu)\frac{\partial u}{\partial x}+\lambda
\frac{\partial v}{\partial y},
\end{equation}
\begin{equation}\label{sig-y}
 \sigma_y =  \lambda \frac{\partial u}{\partial x} +
(\lambda+2\mu)\frac{\partial v}{\partial y},
\end{equation}
\begin{equation}\label{tau-xy}
\tau_{xy}=\mu\left(\frac{\partial v}{\partial x}+\frac{\partial
u}{\partial y}\right).
\end{equation}
Thus, stresses  $\sigma_x$ and $\sigma_y$ are found in $D$ by the
formulas (\ref{sig-x}), (\ref{sig-y}), where
$\mathcal{V}_{1}=\frac{\partial u}{\partial x}$,
$\mathcal{V}_{2}=\frac{\partial v}{\partial y}$. Values of a stress
 $\tau_{xy}$ can not be found without values of the second order partial derivatives
(\ref{W_1-2})  
  in the domain $D$, but the latter can be found with use of Theorem~\ref{syst-1-4} and taking into account a  equivalence
  of the $(u_{x},v_{y})$-problem and the appropriate  (1-4)-problem.

 Indeed, it is well-know (cf., e.g.,\cite[p.,6]{Lu}), that the right-hand side of the equality   (\ref{tau-xy})
 equals
$-\frac{\partial^2 W}{\partial x \partial y}$, consequently,
  $\tau_{xy}=-\frac{\partial^2 W}{\partial x \partial y}$. Therefore,
  a process of finding a stress  $\tau_{xy}$ is reduced to finding the mixed second order partial derivative  $W_{1,1}:=
  \frac{\partial^2 W}{\partial x \partial y}$ in the domain
   $D$. After solving the (1-4)-problem with boundary data
   (\ref{u_1-4}), we rebuild
  in $D$ functions $W_k$, $k=1,2$, in (\ref{W_1-2}).
  We have equalities
\begin{equation}\label{comput-tau-xy}
  \tau_{xy}\equiv \tau_{xy}(x,y)=-W_{1,1}(x,y)=    
  -\int\limits_{(x_{\ast},y_{\ast})}^{(x,y)}\frac{\partial W_{1}}{\partial
  x}\,dx - \int\limits_{(x_{\ast},y_{\ast})}^{(x,y)}\frac{\partial
W_{2}}{\partial y}\,dy  \quad \forall (x,y)\in D,
\end{equation}
where $(x_{\ast},y_{\ast})$ is a fixed point in D, integration means
along any piecewise smooth curve, which joints this point with a
point with variable coordinates $(x, y)\in D$.

Consequently, formulas  (\ref{sig-x}), (\ref{sig-y}),
(\ref{comput-tau-xy}) deliver stresses under boundary conditions
(\ref{u_x-v_y}).

In order to find displacements  $u$ and $v$ we need to obtain their
partial derivatives  $\mathcal{V}_{3}=\frac{\partial u}{\partial
y}$, $\mathcal{V}_{4}=\frac{\partial v}{\partial x}$ in the domain
$D$, so that, displacements are obtained by the formulas
\begin{equation}\label{comput-u-v}
u \equiv
u(x,y)=\int\limits_{(x_{\ast},y_{\ast})}^{(x,y)}\mathcal{V}_{1}\,dx+\mathcal{V}_{3}\,dy,\,\,
v \equiv
v(x,y)=\int\limits_{(x_{\ast},y_{\ast})}^{(x,y)}\mathcal{V}_{4}\,dx+\mathcal{V}_{2}\,dy
\quad \forall (x,y) \in D.
  \end{equation}

Let us find functions  $\mathcal{V}_{3}$ and $\mathcal{V}_{4}$. The
following formulas are fulfilled in $D$ (cf., e.g.,
\cite[p.~6]{Mikhlin},\cite[p.~9]{Lu} with $f_1=f_2 \equiv 0$,
$W:=U$):
\[2\mu u=-\frac{\partial W}{\partial x}+\kappa_0 \left(4p\right),\,
2\mu v=-\frac{\partial W}{\partial y}+\kappa_0 \left(4q\right),
\] where  $W_0+i\,\widetilde{W_0}=4\varphi'(z)$, $\widetilde{W_0}$
is a  harmonic conjugate of  $W_0$ in $D$,
$\varphi(z)=p(x,y)+iq(x,y)$ is an analytic function of the variable $z=x+iy$  in the domain $D_z$. 
These formulas together with the Cauchy--Riemann conditions for
analytic function $\varphi=p+iq$ implies the equalities
\[2\mu \mathcal{V}_{3} 
=-W_{1,1}-\kappa_0 \widetilde{W_0},\,
2 \mu  \mathcal{V}_{4}= 
-W_{1,1}+\kappa_0 \widetilde{W_0}.\] A function $W_{1,1}$ is found
from (\ref{comput-tau-xy}), $\widetilde{W_0}$ is obtained from the
Cauchy--Riemann conditions for analytic function $\varphi=p+iq$.

Finally, formulas (\ref{sig-x}), (\ref{sig-y}),
(\ref{comput-tau-xy}), (\ref{comput-u-v}) define
 the sought-for elastic equilibrium.
\section{Monogenic functions approaches for displacements in 3-D
Elasticity}

 Approaches of monogenic functions with values in non-commutative real algebras
 are using and developing in 3-D
Elasticity. For investigation of the space analog of the equilibrium
Lam\'{e} system  (\ref{Lame}), considered monogenic (in the sense of
some analog of (\ref{usl_K_R})) functions $f \colon \mathcal{A}
\cong \mathbb{R}^{3}\longrightarrow \widetilde{\mathbb{H}}$, where
$\widetilde{\mathbb{H}}\equiv \mathbb{H}:=x_0+x_1 e_1 + x_2 e_2+x_3
e_3$, $x_k\in\mathbb{R}$, $k=\overline{0,4}$,
($e_{k}e_{j}+e_{j}e_{k}=-2\delta_{k,j}$, $k,j=1,2,3$; $e_1 e_2=e_3$)
are real quaternions (cf., e.g.,
\cite{Grig-Trends-14,Grig-Weimar-15,Bock_Gurl}) or
$\widetilde{\mathbb{H}} \equiv \mathcal{A} :=x_0+x_1 e_1 + x_2 e_2
\subset \mathbb{H}$ (cf., e.g.,
\cite{Gur-Greece-14,Gurl-3-D-el+quatern-potenth-14,
g-gen-sol-Lam-syst-mma-15}). In this way, a general regular solution
of the noticed system for a domain with some convex structure is
expressed as a sum of such hypercomplex monogenic functions (for
instance, in \cite{g-gen-sol-Lam-syst-mma-15}: a sum of
(paravector-valued) monogenic, an anti-monogenic and a
$\psi$--hyperholomorphic functions, respectively). So, this
expressions are alternative to the Kolosov-Muskhelishvili formula
for the elastic displacement field and they can be used to solving
boundary value problems for monogenic functions.





\end{document}